
\input amstex
\magnification=1200
\hyphenation{Za-mo-lod-chi-ko-va ko-nech-no-mer-nyh Vig-ne-ra nad-le-zha-shchim
od-no-pa-ra-met-ri-ches-kuyu pa-ra-sta-tis-ti-ki al-geb-ra-i-ches-ko-go
ope-ra-tor-nyh al-geb-ry geo-met-rii pe-re-sta-no-voch-ny-mi urav-ne-ni-yam
ra-bo-te di-na-mi-ches-koe kvad-ra-tich-nye obyk-no-ven-nym
kom-mu-ta-tsi-on-nyh so-ot-no-she-ni-ya-mi di-na-mi-ches-ki na-bo-rom
pe-re-sta-no-voch-nyh}
\font\cyr=wncyr10
\font\cyb=wncyb10
\font\cyi=wncyi10
\font\cyre=wncyr8
\documentstyle{amsppt}
\NoRunningHeads
\NoBlackBoxes
\document
\topmatter\qquad\qquad\qquad\qquad\qquad\qquad\qquad\qquad\qquad\qquad
$\boxed{\boxed{\aligned&\text{\eightpoint "Thalassa Aitheria" Reports}\\
&\text{\eightpoint RCMPI-95/03}\endaligned}}$
\newline
\ \newline
\ \newline
\ \newline
\title\cyr Dinamicheskaya obratnaya zadacha teorii predstavlenii0 i
nekommutativnaya geometriya\endtitle
\author\cyr Denis V. Yurp1ev\endauthor
\address
"Thalassa Aitheria" Research Center for Mathematical Physics and
Informatics
\endaddress
\email denis\@juriev.msk.ru
\endemail
\abstract\nofrills\cyre Dinamicheskaya obratnaya zadacha teorii predstavlenie0,
postanovka kotoroe0 voshodit k klassicheskoe0 rabote E.Vignera ob
opredelyaemosti kommutatsionnyh sootnoshenie0 na kvantovomehaniheskie
velichiny po kvantovym uravneniyam dvizheniya, proillyustrirovana na
prostee0shih primerah.
\endabstract
\endtopmatter
\cyr Teoriya predstavlenie0 imeet po menp1shee0 mere dve storony, odna iz nih
svyazana s {\cyi pryamoe0 zadachee0 teorii predstavlenie0}, zadachee0 opisaniya
predstavlenie0 zadannogo algebraicheskogo obp1ekta. Buduchi v tselom
zavershennoe0 dlya klassicheskih algebraicheskih obp2ektov (algebr i
superalgebr
Li) e1ta problematika vstupila v ocherednoe0 renessans, vyzvannye0 otkrytiem
znachitelp1nogo chisla novyh algebraicheskih struktur v sovremennoe0
kvantovoe0 teorii polya (kvantovye gruppy, algebry Zamolodchikova, operatornye
algebry kvantovoe0 teorii polya i ih raznovidnosti, $W$-algebry i ih
obobshcheniya, kategorii bordizmov i shlee0fy konechnomernyh i
beskonechnomernyh grupp, gomotopicheskie algebry Li, algebry
Batalina--Vilkovyskogo i t.d.). Vtoraya storona svyazana s {\cyi obratnoe0
za\-da\-chee0 teorii predstavlenie0\/}: po imeyushchee0sya sovokupnosti
operatorov
opi\-satp1 "korrektnye0" algebraicheskie0 obp2ekt, imi predstavlyaemye0
(sm.napr. [1]). E1ta zadacha dopuskaet i bolee obshchee rasshirennoe tolkovanie
(sm. [2]).

Razlichnye illyustratsii k obratnoe0 zadache teorii predstavlenie0, a takzhe
ee svyazp1 s nekommutativnoe0 geometriee0 [3], rassmatrivalisp1 v [4,1]. Tselp1
dannoe0 raboty --- vychlenitp1 "dinamicheskuyu" parallelp1 k obratnoe0 zadache,
vos\-hodyashchuyu k klassicheskoe0 rabote E.Vignera [5].

\definition{\cyb Dinamicheskaya obratnaya zadacha teorii predstavlenie0
(zadacha Vignera)} \cyr Dlya zadannoe0 sistemy operatornyh differentsialp1nyh
uravnenie0
$$\frac{\partial}{\partial t}\hat X_i(t)=f_i(\hat X_1(t)\ldots\hat X_n(t))
\qquad (1\le i\le n)\tag *$$
nae0ti algebraicheskie0 obp2ekt $\frak A$, porozhdennye0 generatorami $e_i$ i
odnoparametricheskuyu gruppu $\gamma_t$ avtomorfizmov $\frak A$ takie, chto
pri lyubom (nadlezhashchim obrazom opredelennom) operatornom predstavlenii
$T$ algebraicheskogo obp2ekta $\frak A$ gruppa $\gamma_t$ opredelyaet\-sya
sistemoe0 uravnenie0 $(*)$, kolp1 skoro $T(e_i)=\hat X_i(0)$.
\enddefinition

\cyr Algebraicheskie0 obp2ekt $\frak A$ mozhet bytp1 assotsiativnoe0 algebroe0
polinomialp1nogo rosta, hotya i drugie varianty vozmozhny (naprimer,
parastatistiki [5,6] ili sootnosheniya v izotopicheskih parah [7,8]).
V dalp1\-nee0\-shem, odnako, ukazannye vozmozhnosti ne budut rassmatrivatp1sya.

\cyr Prakticheski naibolp1shie0 interes predstavlyayut beskonechnye sistemy
operatornyh uravnenie0, otvechayushchie nelinee0nym operatornym uravneniyam
v chastnyh proizvodnyh kvantovoe0 teorii polya (napr. uravneniyam sae0n--Gordon
ili Liuvillya, rassmatrivavshimsya s ukazannoe0 tochki zreniya v kachestve
harakternyh primerov leningradskoe0 shkoloe0 akad. L.D.Faddeeva v ramkah
kvantovogo metoda obratnoe0 zadachi, gruppoe0 iz Vysshee0 Normalp1noe0 Shkoly
v Parizhe pod rukovodstvom prof.Zh.-L.Zher\-ve1, kollektivom otdela kvantovoe0
teorii polya v Institute Steklova v Moskve). Odnako, s metodologicheskoe0
tochki zreniya tselesoobrazno rassmatrivatp1 i konechnye sistemy operatornyh
uravnenie0, otvechayushchie obyknovennym differentsialp1nym uravneniyam
kvantovoe0 mehaniki.

\definition{\cyb Opredelenie 1} \cyr {\cyi Dinamicheskim osvobozhdeniem skrytyh
simmetrie0\/} na\-zy\-va\-et\-sya sopostavlenie zadannoe0 sisteme operatornyh
dif\-fe\-ren\-tsi\-alp1\-nyh uravnenie0 $(*)$ pary $(\frak A,\gamma_t)$, gde
$\frak A$ --- assotsiativnaya algebra, porozhdennaya generatorami $e_1,\ldots
e_n$ tak, chto otobrazhenie simmetrizatsii Vee0lya $w:S^{\cdot}(H)\mapsto\frak
A$ ($H=\operatorname{span}(e_1,\ldots e_n)$) yavlyaet\-sya izomorfizmom
linee0nyh prostranstv, a $\gamma_t$ --- odnoparametricheskaya gruppa
avtomorfizmov algebry $\frak A$ takaya, chto pri lyubom operatornom
predstavlenii $T$ assotsiativnoe0 algebry $\frak A$ gruppa $\gamma_t$
opredelyaet\-sya sistemoe0 $(*)$, kolp1 skoro $T(e_i)=\hat X_i(0)$.
Dinamicheskoe osvobozhdenie simmetrie0 nazyvaet\-sya {\cyi gamilp1tonovym\/},
esli sushchestvuet takoe0 e1lement $h$ algebry $\frak A$, chto sistema $(*)$
privodit\-sya k vidu $\frac{\partial}{\partial t}\hat X_i=[\hat H,\hat X_i]$
($\hat H=T(h)$).
\enddefinition

\cyr S algebrogeometricheskoe0 tochki zreniya dinamicheskoe osvobozhdenie
\linebreak skrytyh simmetrie0 zaklyuchaet\-sya v interpretatsii zadannoe0
sistemy
operatornyh differentsialp1nyh uravnenie0 kak potoka (dinamicheskoe0 sistemy)
na nekommutativnom algebraicheskom mnogoobrazii [3]. Ne\-ob\-ho\-di\-mo
otmetitp1,
chto v otlichie ot kommutativnogo sluchaya dannaya interpretatsiya ne
yavlyaet\-sya trivialp1noe0 i samoochevidnoe0. Zametim takzhe, chto
dinamicheskie
sistemy na nekommutativnyh mnogoobraziyah voznikayut estestvennym obrazom i v
ramkah nekommutativnoe0 differentsialp1noe0 geometrii Alana Konna [9].

\cyr Otmetim, chto interpretatsiya kvantovoe0 konformnoe0 teorii polya,
tesno svyazannoe0 s upomyanutymi vyshe uravneniyami kvantovoe0 teorii polya,
kak beskonechnomernoe0 nekommutativnoe0 geometrii byla dana avtorom v rabote
[10]. Odnako, podcherknem, chto v konkretnyh kvantovopolevyh zadachah
zachastuyu ispolp1zuyut\-sya assotsiativnye algebry $\frak A$ s fiksirovannym
naborom e1lementov $e_i$ takih, chto $T(e_i)=\hat X_i(0)$, no ne porozhdaemyh
imi. Kak pravilo, v kachestve algebry $\frak A$ udobno rassmatrivatp1
algebry Vee0lya ot koe1ffitsientov svobodnyh polee0 (t.n. {\cyi predstavlenie
svobodnyh polee0}), universalp1nye obertyvayushchie algebry algebr Katsa--Mudi
ili nekotorye ih $q$--deformatsii. Tem samym podobnye konstruktsii
yavlyayut\-sya estestvennymi alp1ternativami dinamicheskomu osvobozhdeniyu
skrytyh simmetrie0 pri reshenii dinamicheskoe0 obratnoe0 zadachi teorii
predstavlenie0.

\cyr Rassmotrim prostee0shie primery dinamicheskogo osvobozhdeniya skrytyh
simmetrie0.

\example{\cyb Primer 1 ({\cyr $n=2$, $f_1$ i $f_2$ --- linee0nye funktsii})}
\cyr V e1tom sluchae sistema $(*)$ imeet vid:
$$\left\{\aligned
\frac{\partial}{\partial t}\hat X&=a\hat X+b\hat Y\\
\frac{\partial}{\partial t}\hat Y&=c\hat X+d\hat Y
\endaligned\right.$$
Bez ogranicheniya obshchnosti mozhno schitatp1, chto matritsa $A=\left(
\matrix a & b \\ c & d \endmatrix\right)$ privedena k zhordanovoe0 normalp1noe0
forme. Pustp1 matritsa $A$ diagonalp1na: $A=\left(\matrix \lambda & 0 \\
0 & \mu \endmatrix\right)$. Sistema $(*)$ zapisyvaet\-sya kak
$$\left\{\aligned
\frac{\partial}{\partial t}\hat X&=\lambda\hat X\\
\frac{\partial}{\partial t}\hat Y&=\mu\hat Y
\endaligned\right.$$
Budem iskatp1 sootnosheniya na $\hat X$ i $\hat Y$ v vide $[\hat X,\hat Y]=
\hat f(\hat X,\hat Y)$, gde $\hat f(\hat X,\hat Y)=\sum_{i\ge 0}\sum_{j\ge 0}
a_{ij}\hat X^i\hat Y^j$. Uslovie sohraneniya dinamikoe0 ukazannyh
so\-ot\-no\-she\-nie0 imeet vid
$$(\lambda+\mu)\sum_{i\ge o}\sum_{j\ge 0}a_{ij}\hat X^i\hat Y^j=
\sum_{i\ge 0}\sum_{j\ge 0}a_{ij}(i\lambda+j\mu)\hat X^i\hat Y^j$$
ili
$$\sum_{i\ge 0}\sum_{j\ge 0}a_{ij}(\lambda(i-1)+\mu(j-1))\hat X^i\hat Y^j=0.$$
Kak sledstvie, koe1ffitsienty $a_{ij}$ otlichny ot nulya tolp1ko dlya teh
$i$ i $j$, dlya kotoryh vyrazhenie $\lambda(i-1)+\mu(j-1)$ obrashchaet\-sya
v nulp1.

Vozmozhny pyatp1 sluchaev:
\roster
\item"1)" Esli $\lambda$ i $\mu$ nesoizmerimy, ili $\lambda$ i $\mu$
soizmerimy,
odnogo znaka i $\lambda\ne\mu$, to $[\hat X,\hat Y]=\alpha\hat X\hat Y$ ili,
inymi slovami, $\hat X\hat Y=q\hat Y\hat X$ (sootnosheniya dvymernoe0
kvantovoe0
ploskosti [3]).
\item"2)" Esli $\lambda=\mu$, to $[\hat X,\hat Y]=\alpha\hat X^2+\beta\hat X
\hat Y+\gamma\hat Y^2$ (algebra s kvadratichnymi sootnosheniyami, porozhdennaya
dvumya obrazuyushchimi [3,11,12]).
\item"3)" Esli $\lambda=0$, $\mu\ne 0$, to $[\hat X,\hat Y]=P(\hat X)\hat Y$,
ili, inymi slovami, $\hat Y\hat X=R(\hat X)\hat Y$ (algebry podobnogo tipa
rassmatrivalisp1 v [11,12,13]).
\item"4)" Esli $\lambda+\mu=0$, to $\hat X,\hat Y]=\sum_{k\ge 0}a_k\hat X^k
\hat Y^k$. E1tot sluchae0 otvechaet klassicheskoe0 zadache Vignera [5] ob
odnoznachnosti vosstanovleniya kvantovyh kommutatsionnyh sootnoshenie0 po
uravneniyam kvantovoe0 dinamiki. V chastnosti, dannye0 sluchae0 vklyuchaet
t.n. $q$--kom\-mu\-ta\-tsi\-on\-nye sootnosheniya
$\hat X\hat Y-q\hat Y\hat X=1$ (sm. napr. [14]).
\item"5)" Esli $\lambda$ i $\mu$ soizmerimy, $\lambda=l\nu$, $\mu=m\nu$,
$l<0$, $m>0$, $l+m\ne 0$, to $[\hat X,\hat Y]=\sum_{k\ge 0}a_k\hat X^{1+km}
\hat Y^{1-kl}$.
\endroster
Otmetim, chto dinamicheskoe osvobozhdenie skrytyh simmetrie0 yavlyaet\-sya
gamilp1tonovym lishp1 v sluchae kanonicheskih kommutatsionnyh sootnoshenie0,
yavlyayushchihsya podsluchaem sluchaya 4, i linee0nyh kommutatsionnyh
sootnoshenie0, yavlyayushchihsya podsluchaem sluchaya 3 s $R(x)\equiv 1$.

Pustp1 teperp1 matritsa $A$ imeet vid netrivialp1noe0 zhordanovoe0 kletki
$\left(\matrix \lambda & 1 \\ 0 & \lambda\endmatrix\right)$, izuchim, kogda
sushchestvuet dinamicheskoe osvobozhdenie skrytyh simmetrie0 s algebroe0
$\frak A$, zadavaemoe0 linee0nymi ili kvadratichnymi kommutatsionnymi
sootnosheniyami na generatory. Sistema $(*)$ imeet vid:
$$\left\{\aligned\frac{\partial}{\partial t}\hat X&=\lambda\hat X\\
\frac{\partial}{\partial t}\hat Y&=\lambda\hat Y+\hat X\endaligned\right.$$
\roster
\item"1)" Pri $\lambda\ne 0$ ne sushchestvuet algebry $\frak A$ dinamicheski
osvobozhdennyh skrytyh simmetrie0 s linee0nymi kommutatsionnymi sootnosheniyami
na generatory.
\item"2)" Kvadratichnye kommutatsionnye sootnosheniya v algebre $\frak A$
dinamicheski osvobozhdennyh skrytyh simmetrie0 predstavlyayut\-sya v vide
$[\hat X,\hat Y]=a\hat Y^2$. Sootvet\-stvuyushchee dinamicheskoe osvobozhdenie
skrytyh simmetrie0 ne yavlyaet\-sya gamilp1tonovym.
\endroster
\endexample

\example{\cyb Primer 2 ({\cyr $n=2$, $f_1$ i $f_2$ --- odnorodnye kvadratichnye
funktsii})}\cyr Sistema $(*)$ imeet vid:
$$\left\{\aligned
\frac{\partial}{\partial t}\hat X&=a\hat X^2+b\tfrac{\hat X\hat Y+\hat Y\hat
X}2+
c\hat Y^2\\
\frac{\partial}{\partial t}\hat Y&=d\hat X^2+e\tfrac{\hat X\hat Y+\hat Y\hat
Y}2+
f\hat Y^2
\endaligned\right.\tag {**}$$
Izuchim, kogda dannaya sistema operatornyh uravnenie0 dopuskaet dinamicheskoe
osvobozhdenie skrytyh simmetrie0 s algebroe0 $\frak A$, opredelyayushchie
kommutatsionnye sootnosheniya kotoroe0 odnorodnye linee0nye ili kvadratichnye
funktsii, a takzhe kogda dinamicheskoe osvobozhdenie skrytyh simmetrie0
yavlyaet\-sya gamilp1tonovym.
\roster
\item"1)" Sistema operatornyh uravnenie0 $(**)$ dopuskaet dinamicheskoe
osvobozhdenie skrytyh simmetrie0 s algebroe0 $\frak A$, opredelyayushchie
kommutatsionnye sootnosheniya kotoroe0 linee0ny, togda i tolp1ko togda, kogda
$$\left\{\aligned
ac+bf+f^2&=0\\
a^2+ae+df&=0
\endaligned\right.$$
Pri e1tom $[\hat X,\hat Y]=\alpha\hat X+\beta\hat Y$, i imeet mesto proportsiya
$\alpha\!:\!\beta=-a\!:\!f$. Dinamicheskoe osvobozhdenie skrytyh simmetrie0
gamilp1tonovo, $\hat H=p\hat X^2+q\tfrac{\hat X\hat Y+\hat Y\hat X}2+r\hat Y^2$
i imeet mesto proportsiya $p\!:\!q\!:\!r=df\!:\!af\!:\!ac$.
\item"2)" Sistema operatornyh uravnenie0 $(**)$ dopuskaet dinamicheskoe
osvobozhdenie skrytyh simmetrie0 s algebroe0 $\frak A$, opredelyayushchie
kommutatsionnye sootnosheniya kotoroe0 kvadratichny, togda i tolp1ko togda,
kogda $a\!:\!d=b\!:\!e=c\!:\!f$. Pri e1tom $[\hat X,\hat Y]=\alpha\hat
X^2+\beta\frac{\hat
X\hat Y+\hat Y\hat X}2+\gamma\hat Y^2$, i imeyut mesto proportsii
$a\!:\!d\!:\!\alpha=b\!:\!e\!:\!\beta=c\!:\!f\!:\!\gamma$. Dinamicheskoe
osvobozhdenie skrytyh
simmetrie0 gamilp1tonovo, $\hat H=p\hat X+q\hat Y$ i
$p\!:\!q=-a\!:\!d=-b\!:\!e=-c\!:\!f$.
\endroster
Otmetim, chto sistema $(*)$ dopuskaet dinamicheskoe osvobozhdenie skrytyh
simmetrie0 kak s linee0nymi, tak i s kvadratichnymi opredelyayushchimi
kommutatsionnymi sootnosheniyami v algebre $\frak A$, esli $b=-(a+c)$, $d=a$,
$e=b$, $f=c$. Linee0nye0 sluchae0 zadaet\-sya sootnosheniyami $[\hat X,\hat
Y]=a\hat X+b\hat Y$ i gamilp1tonianom $\hat H=\frac12(\hat X+\hat Y)^2$,
kvadratichnye0 --- sootnosheniyami $[\hat X,\hat Y]=a\hat X^2-(a+c)\tfrac{\hat
X\hat Y+\hat Y\hat X}2+c\hat Y^2$ i gamilp1tonianom $\hat H=\hat X-\hat Y$.
\endexample

\cyr Interesnye vozmozhnosti voznikayut v sluchae sistemy iz treh operatornyh
differentsialp1nyh uravnenie0 s linee0noe0 ili kvadratichnoe0 pravoe0
chastp1yu,
odnako, on ne soderzhit ideologicheski nichego printsipialp1no
novogo, esli ogranichivatp1sya algebrami $\frak A$ dinamicheski osvobozhdennyh
skrytyh simmetrie0 s linee0nymi ili kvadratichnymi perestanovochnymi
sootnosheniyami na generatory.

Otmetim, chto i pri lyubom konechnom $n$ issledovanie vozmozhnosti
dinamicheskogo
osvobozhdeniya skrytyh simmetrie0, opisyvaemyh linee0nymi ili kvadratichnymi
perestanovochnymi sootnosheniyami, v sistemah $n$ operatornyh
differentsialp1nyh
uravnenie0 s linee0noe0 ili kvadratichnoe0 pravoe0 chastp1yu svodit\-sya k
voprosu o sushchestvovanii netrivialp1nyh re\-she\-nie0 pereopredelennyh sistem
linee0nyh uravnenie0 s posleduyushchee0 pro\-ver\-koe0 tozhdestv Yakobi i ne
predstavlyaet teoreticheskih ili vychislitelp1nyh trudnostee0 v otlichie ot
sluchaya silp1no nelinee0nyh perestanovochnyh sootnoshenie0 (sm. [12]).

\head\cyb Spisok literatury\endhead

\roster\eightpoint
\item"[1]" Juriev D., Topics in hidden symmetries. E--print (LANL Electronic
Archive on Theor. High Energy Phys.): {\it hep-th/9405050} (1994).
\item"[2]" {\cyre Yurp1ev D.V., Harakteristiki par operatorov, gibridy Li,
skobki Puassona i nelinee0naya geometricheskaya algebra}. Report RCMPI-95/02
(1995) and e--print (SISSA Electronic Archive on Funct. Anal.):
{\it funct-an/9411007} (1994).
\item"[3]" Manin Yu.I., {\it Topics in noncommutative geometry},
Princeton Univ. Press, Princeton, NJ, 1991.
\item"[4]" Juriev D., Setting hidden symmetries free by the noncommutative
Veronese mapping, J. Math. Phys. 35 (1994) 5021-5024.
\item"[5]" Wigner E.P., Phys. Rev. 77 (1950) 711-712.
\item"[6]" Ohnuki Y., Kamefuchi S., {\it Quantum field theory and
parastatistics}, Springer, 1982.
\item"[7]" Juriev D., Topics in isotopic pairs and their representations,
Theor. Matem. Fiz. 00 (1995) 00-00 and e--print (Texas Electronic Archive
on Math. Phys.): {\it mp\_arc/94-267} (1994).
\item"[8]" Juriev D., Classical and quantum dynamics of noncanonically
coupled oscillators and Lie superalgebras, Russian J. Math. Phys. 00 (1995)
00-00 and e--print (SISSA Electronic Archive on Funct. Anal.): {\it
funct-an/9409003} (1994).
\item"[9]" Connes A., Non--commutative differential geometry. Publ. Math.
IHES 62 (1986).
\item"[10]" {\cyre Yurp1ev D.V., Kvantovaya konformnaya teoriya polya kak
beskonechnomernaya nekommutativnaya geometriya. Uspekhi Matem. Nauk 46(3)
(1991) 115-138.}
\item"[11]" {\cyre Karasev M.V., Maslov V.P., Nazae0kinskie0 V.E., Algebry s
obshchimi perestanovochnymi sootnosheniyami. {\rm I}, {\rm II}. Sovrem. probl.
matem. Novee0sh. dostizheniya. T.13. M., VINITI, 1979}.
\item"[12]" {\cyre Karasev M.V., Maslov V.P., Nelinee0nye skobki Puassona.
Geometriya i kvantovanie. M., Nauka, 1991}.
\item"[13]" {\cyre Latushkin Yu.D., Stepin A.M., Uspekhi Matem. Nauk
46 (1991) 85}.
\item"[14]" {\cyre Kulish P.P., Kontraktsiya kvantovyh algebr i
$q$--ostsillyatory. Teor. Matem. Fiz. 86(1) (1991) 157-160}.
\endroster
\enddocument